\documentclass{amsart}

\usepackage[all]{xy}

\usepackage{pictexwd, epsfig, amsthm, amssymb, amsmath, graphicx, psfrag,amsxtra, latexsym}
\usepackage{enumerate}

\newcommand{\mD}{\mathbb D}
\newcommand{\mZ}{\mathbb Z}
\newcommand{\mR}{\mathbb R}
\newcommand{\mC}{\mathbb C}

\DeclareMathOperator{\im}{im}

\theoremstyle{plain}
\newtheorem{Thm}{Theorem}

\newtheorem*{lemma}{Lemma}

\theoremstyle{definition}

\theoremstyle{remark}
\newtheorem*{Prf}{Proof}

\begin{document}

\title[A note on strong Jordan separation.]
      {A note on strong Jordan separation.}
\author{Jean-Fran\c{c}ois\ Lafont}
\address{Department of Mathematics\\
         Ohio State University\\
         Columbus, OH  43210}
\email[Jean-Fran\c{c}ois\ Lafont]{jlafont@math.ohio-state.edu}

\begin{abstract}
We establish a strengthening of Jordan separation, to the setting of
maps $f: X \rightarrow S^{n+1}$, where $X$ is not necessarily a manifold,
and $f$ is not necessarily injective. 
\end{abstract}

\maketitle

\section{Introduction}

In a previous paper \cite{L}, we established a result which we termed 
{\it strong Jordan separation}. This was a version of Jordan separation which
applied to maps $f:S^n\rightarrow S^{n+1}$ which aren't assumed to be injective.
Under some mild hypothesis, one could nevertheless ensure that the image 
separated $S^{n+1}$, and that any continuous extension 
$F: \mD^{n+1}\rightarrow S^{n+1}$ surjects onto one of the connected components
of $S^{n+1}-f(S^n)$. Recently, Iwaniec and Onninen \cite{IO} found applications of this result
in the field of quasi-conformal hyperelasticity. In this short note, we extend our 
result from \cite{L} to the broadest possible setting, by establishing the following two
results:

\begin{Thm} Let $X$ be a compact topological space, $f:X\rightarrow S^{n+1}$ a continuous map,
and $U\subset X$ an open subset homeomorphic
to an open $n$-disk $\mD^n_\circ$.
Assume that 
\begin{itemize}
\item the map $f: X\rightarrow S^{n+1}$ contains $U$ in its set
of injectivity (i.e. $U\subset Inj(f):= \{x\in X \hskip 2pt | \hskip 2pt f^{-1}(f(x))=x\}$), and
\item the map $\check{H} ^n(X; \mZ_2)\rightarrow \check{H} ^n(X-U; \mZ_2)$ on
 \v{C}ech cohomology groups  induced by the inclusion $X-U\hookrightarrow X$ 
 has a non-trivial kernel.
\end{itemize}
Then $f(X)$ separates $S^{n+1}$ into at least two connected components. Furthermore, 
there are precisely two connected components $V_1,V_2$ of $S^{n+1}-f(X)$ having
the property that their closure $\bar V_i$ intersects $f(U)$. In fact, for these two 
connected components, we have containments $f(U) \subset \bar V_i$.
\end{Thm}

In the previous theorem, one should think of the sets 
$V_1,V_2$ as corresponding locally to the two ``sides'' of $f(U)\cong 
\mD^n_\circ$ in the ambient $S^{n+1}$.

\begin{Thm}
Under the hypotheses of the previous theorem, let us further assume that $X$ is a closed
subspace of an ambient topological space $\hat X$. 
Define two subgroups of $H_n(X; \mZ_2)$ by:
\begin{itemize}
\item $K=\ker \big(H_n(X; \mZ_2)\rightarrow 
H_n(\hat X ; \mZ_2)\big)$, and 
\item $J= \im \big(H_n(X-U; \mZ_2)\rightarrow H_n(X ; \mZ_2)\big)$.
\end{itemize} 
where both maps are induced by the corresponding inclusions of spaces. 
If $K\not \subseteq J$, then we have that for 
any continuous extension $F:\hat X\rightarrow S^{n+1}$, $F$ surjects onto one of the two
open components $V_i$.
\end{Thm}

We will prove both these theorems in Section 2 of our paper.  Observe that,
other than the condition requiring $U$ to lie in $Inj(f)$, the hypotheses in both theorems
are {\it internal},
in the sense that they are statements purely about the (homology of the) 
spaces $U, X, \hat X$, and do not involve the map $f$. In Section 3,
we will provide examples showing that the hypotheses of the theorems are necessary. We will
also exhibit examples of triples $(\hat X, X, U)$ satisfying the homological hypotheses of both theorems.

\vskip 10pt

\centerline{\bf Acknowledgments}

\vskip 5pt

Tadeusz Iwaniec asked whether the author's arguments for the strong Jordan separation 
theorem in \cite{L} could be extended to the case where the source space is a manifold other than the 
sphere. The present paper was motivated by, and answers, Iwaniec's question.

The author was partially supported by the NSF, under grant DMS-0606002, and by a
Sloan Research Fellowship.

\section{Proofs}

Throughout this section, all homology and cohomology groups are understood to have
coefficients in $\mZ_2$. The proofs of both Theorems follow closely the proof of the strong
Jordan separation in \cite[Section 2]{L}.

\begin{Prf}[Theorem 1]
We start by recalling that Alexander duality provides us with an isomorphism:
$$\tilde H_0(S^{n+1} - f(X)) \cong \check H ^n(f(X))$$
hence to show that $f(X)$ separates $S^{n+1}$ it is sufficient to show that 
$ \check H^n(f(X)) \neq 0$. Now consider the 
decomposition of $X$ into two open sets, one obtained by shrinking $U$ 
slightly, the other obtained by enlarging $X-U$ slightly. The intersection of these
two open sets is an open subset homeomorphic to $S^{n-1}\times \mR$. Corresponding
to this decomposition, we have an image decomposition of $f(X)$ into two open sets
(recall that $U\subset Inj(f)$). In particular, we can apply the Mayer-Vietoris sequence
in \v Cech cohomology to compute the cohomology of $f(X)$:
$$ \ldots \rightarrow \check H^{n-1}(S^{n-1})\rightarrow \check H^n(f(X)) \rightarrow 
\check H^n(f(U))\oplus \check H^n(f(X-U)) \rightarrow \check H^n(S^{n-1}) \rightarrow \ldots$$
Putting in the known terms into the sequence above, we obtain:
\begin{equation}
\ldots \rightarrow \mZ_2 \rightarrow \check H^n(f(X)) \rightarrow \check H^n(f(X-U)) \rightarrow 0
\end{equation}
so to show that $ \check H^n(f(X)) \neq 0$, it is sufficient to show that the $\mZ_2 \cong 
\check H^{n-1}(S^{n-1})$ injects into $ \check H^n(f(X))$. In order to show this, we compare
the Mayer-Vietoris sequence above with the corresponding Mayer-Vietoris sequence for
the decomposition of $X$. This gives us the following commutative diagram:
$$\xymatrix{
\ldots \ar[r]  &  \mZ_2 \ar[r]^-{\phi} \ar[d]^{f^*}_{\cong} & \check H^n(f(X)) \ar[r] \ar[d]^{f^*} & 
\check H^n(f(X-U)) \ar[r] \ar[d]^{f^*} & 0\\
\ldots \ar[r] &  \mZ_2 \ar[r] ^-\psi & \check H^n(X) \ar[r] & \check H^n(X-U)
\ar[r] & 0 \\
}$$
Since we assumed that the map $\check H^n(X) \rightarrow \check H^n(X-U)$ has
non-trivial kernel, we see that $\psi$ is injective. This forces the composite 
$\psi \circ f^* = f^*\circ \phi$ to be injective, and hence the map $\phi$ to be injective,
as desired. This allows us to conclude that 
$0 \neq \check H ^n(f(X)) \cong \tilde H_0(S^{n+1} - f(X))$, and we see that 
$f(X)$ separates $S^{n+1}$. 

Now, the remainder of the proof is virtually identical to that given in \cite{L}.
We will briefly sketch out the arguments, referring the reader to \cite[Section 2.1]{L}
for more details. 

To see that there are precisely two connected components $V_1,V_2$ whose closure 
intersects $f(U)$, one considers
the inclusion $S^{n+1}-f(X) \hookrightarrow S^{n+1}-f(X-U)$. The latter space is obtained
from the former by ``adding in $f(U)$'', i.e. we have identifications 
$S^{n+1}-f(X-U) = \big( S^{n+1} - f(X)\big) \cup f(U)$. By applying Alexander duality to the
exact sequence in (1), we obtain:
$$0\rightarrow \mZ_2 \rightarrow \tilde H_0(S^{n+1}-f(X))\rightarrow  
\tilde H_0\big(\big( S^{n+1} - f(X)\big) \cup f(U)\big)\rightarrow 0$$
This is a homological version of the statement ``there are precisely two components 
of $S^{n+1}-f(X)$ which are incident to $f(U)$''. 

In order to obtain the statement we desire,
we make use of the following elementary result from point set topology: 
if $\{V_i\}$ is a collection of pairwise disjoint open sets in $\mR^{n+1}$, and $Z$ is a connected
set which intersects the closure of each $V_i$, then $Z \cup \big( \bigcup V_i\big)$ is connected.
Now apply this to the situation where the $\{V_i\}$ are the connected components of $S^{n+1}-
f(X)$ whose closure intersects $f(U)$, and $Z=f(U)$. For any distinct pair of connected components
of $S^{n+1}-f(X)$ whose closure intersects $f(U)$, we will get a corresponding element in 
$\tilde H_0(S^{n+1}-f(X))$ which maps to zero in $H_0(S^{n+1}-f(X-U))$. In other words, the 
rank of the kernel is one less than the number of connected components of $S^{n+1}-f(X)$ 
whose closure intersects $f(U)$. Since we know that the kernel has rank one, we conclude that 
there are precisely two connected components $V_1,V_2$ of $S^{n+1}-f(X)$ whose closure 
intersects $f(U)$. 

Finally, to see that $f(U) \subset \bar V_i$, take $p\in f(U)$ arbitrary, consider a shrinking sequence
of open metric balls $\{U_i\}$ centered at $f^{-1}(p)$. We can apply the same argument as
in the previous paragraph, but replacing $f(U)$ by $f(U_i)$. Observe that $f(X-U)\subset f(X-U_i)$ 
induces an isomorphism on all the \v Cech cohomology groups,
and so by Alexander duality, the homology groups of the complements are also unchanged. 
This forces {\it the same
two components} $V_1,V_2$ to intersect {\it every $f(U_i)$}. Since the sets $f(U_i)$ are shrinking
down to $\{p\}$, this immediately gives us that 
$p$ lies in the closure of both $V_i$, completing the proof of Theorem 1.
\end{Prf}

\begin{Prf}[Theorem 2] Before starting with the proof of the theorem, let us briefly discuss some
general background material. For $p\notin f(X)$, we will consider the homomorphism
$f_*:H_n(X)\rightarrow H_n(\mR^{n+1} - p)\cong \mZ_2$. Note that, since all the groups 
$H_n(\mR^{n+1} - p)$ are isomorphic to $\mZ_2$, we see that for any $p,q \notin f(X)$, 
there are canonical identifications between the groups $H_n(\mR^{n+1} - p)$ and 
$H_n(\mR^{n+1} - q)$. In particular, it makes sense to talk about elements being ``the same'' 
or ``different'' in the groups $H_n(\mR^{n+1} - p)$ and $H_n(\mR^{n+1} - q)$. Finally,
let us fix some notation. Recall that $U\subset X$ is an open set
homeomorphic to an open disc $\mD ^n_\circ$, which we identify with the unit disk in $\mR^n$.
Fixing this identification, we now denote by $U(r)$ ($r<1$) the subset of $U$ which 
corresponds to the open disk of radius $r$. We will use $O$ to denote the point in $U$
which corresponds to the origin.

Let us now argue by way of contradiction: assume that there exists a continuous 
extension $F: \hat X \rightarrow S^{n+1}$ and points 
$z_i\in V_i$ with $z_i\notin F(\hat X)$. By hypothesis, $K\not \subseteq J$, so 
there exists a homology class $\alpha \in H_n(X)$ having the following two properties:
\begin{enumerate}
\item $\alpha \in \ker( H_n(X) \rightarrow H_n(\hat X))$, and
\item $\alpha \notin \im (H_n(X-U) \rightarrow H_n(X)).$
\end{enumerate}
Let us consider the image of the class $\alpha$ in each of the homology groups 
$H_n(\mR^{n+1} - z_i) \cong \mZ_2$. To compute this, observe that since 
$z_i\notin F(\hat X)$, we can
factor the map $f_*$ as the composition $F_*\circ i_*$, where $i: X\rightarrow \hat X$
is the inclusion. This forces the containment $\ker (i_*)\subseteq \ker (f_*)$, which
combined with property (1) in the choice of $\alpha$, ensures that
$\alpha$ maps to zero in each of the homology groups $H_n(\mR^{n+1} - z_i) \cong \mZ_2$. 

So in order to obtain a contradiction, it is sufficient to show that the class $\alpha$ maps to {\it distinct} 
elements in each of the two homology groups $H_n(\mR^{n+1} - z_i) \cong \mZ_2$. This is
considerably harder; we will proceed in several steps. We will first 
replace the map $f$ by
another map $g$, obtained by locally perturbing $f$ on the interior of $U$. This new map $g$
will have the following properties:
\begin{itemize}
\item $g$ is tame on the interior of $U(1/2)$, 
\item $g\equiv f$ on the complement of $U(3/4)$,
\item $U \subset Inj(g)$,
\item $f$ is homotopic to $g$ in the complement of the points $z_1,z_2$.
\end{itemize}
The construction of $g$ can be done by appealing to the important {\it codimension one 
taming theorem} of
Ancel-Cannon \cite{AC} (when $n\geq 4$), Ancel \cite{A} (when $n=3$), Bing \cite{B} (when
$n=2$) and Schoenflies (when $n=1$). We refer the reader to \cite[pg. 689]{L} for details on
how to accomplish this perturbation.

Since $f\simeq g$ in the complement of the $z_i$, it is sufficient to show that $g_*$
maps $\alpha$ to distinct elements in $H_n(\mR^{n+1} - z_i)$. Now for the map $g$, we 
can again apply Theorem 1, and see that there are exactly two connected components 
$V_1^\prime, V_2^\prime$ of $S^{n+1}-g(X)$ whose closure intersects $g(U)$. 
Recall that $O\in U$ is the point corresponding to the origin
under the homeomorphism identifying $\mD^n_\circ$ with $U$. We know that there 
exists a sequence of points $x_i \in V_1^\prime$, $y_i \in V_2^\prime$, with the property
that $\lim \{x_i\}=\lim\{y_i\}=g(O)$. We claim that, for each $i$, the 
$g_*(\alpha)$ form distinct elements in the two groups
$H_n(\mR^{n+1}-x_i)$ and $H_n(\mR^{n+1}-y_i)$. 

To see this, we first note that, if the points $p,q$
are chosen in the same connected component of $\mR^{n+1}-g(X)$, then the two homomorphisms 
$g_*: H_n(X)\rightarrow H_n(\mR^{n+1} - p)$ and $g_*: H_n(X)\rightarrow H_n(\mR^{n+1} - q)$
are identical. In particular, it is enough to show that the images of $g_*(\alpha)$ are distinct for
a specific index $i$. But recall that the map $g$ is {\it tame} on $U(1/2)$, and hence there
exists a global homeomorphism $\phi: \mR^{n+1}\rightarrow \mR^{n+1}$ with the property that
$\phi\circ g$ maps $U(1/2)$ into the standard $\mR^n \times \{0\} \subset \mR^{n+1}$. Furthermore,
for $i$ sufficiently large, $\phi(x_i)$ (respectively $\phi(y_i)$) are points which are locally 
immediately above (respectively below) the hyperplane $\mR^n\times \{0\}$. It is of course
sufficient to show that $\phi_*(g_*(\alpha))$ represents distinct elements in the two homology 
groups $H_n(\mR^{n+1}-\phi(x_i))$ and $H_n(\mR^{n+1}-\phi(y_i))$.

\vskip 10pt

Now consider a cycle representing the homology class $\alpha \in H_n(X)$; this is a formal linear
sum $\sum \sigma_k$ of finitely many maps $\sigma_k: \Delta ^n \rightarrow X$. Now the cycle 
we are interested in is the image cycle $\phi_*(g_*(\alpha))$, which is represented by the formal
linear sum $\sum \tau_k$, where each $\tau_k= \phi\circ g\circ \sigma_k$. In particular, there
exists a point $p \in \phi(g(U(1/2)))$ with the property that $p$ lies solely in the image of the 
{\it interior} of the finitely many simplices, and $p$ is a non-singular value of each of the
maps  $\tau_k$. Now note that
we can join $\phi(x_i)$ to $\phi(y_i)$ by a PL-curve $\eta$ which intersects $(\phi\circ g)(X)$ in a 
single transverse intersection at the point $p$.
There are now two possibilities: either (A) the homology class $\phi_*(g_*(\alpha))$ represents
distinct elements in the two homology groups $H_n(\mR^{n+1}-\phi(x_i))$ and 
$H_n(\mR^{n+1}-\phi(y_i))$, or (B) from intersection theory, we have that the number of
intersection points of $\eta$ with the maps $\tau_k$ is even. 

Let us now argue that possibility (B) does not occur. By way of contradiction, if this was the
case, one could subdivide the finitely many simplices which intersect $\eta$, obtaining a 
new cycle $\sum \tau_k^\prime$ having the additional property that all the maps $\tau_k^\prime$
whose image passes through $p$ coincide with a fixed map $\tau: \Delta ^n\rightarrow 
\mR^{n+1}$. At the cost of further subdividing, we may moreover assume that $\tau$ lies
entirely within the image of $g(U)$. Now, since there are an {\it even} number of copies 
of $\tau$ in the cycle
$\sum \tau_k^\prime$, and since we are working with $\mZ_2$-coefficients, we can remove
all occurences of this singular simplex from the cycle, resulting in a new cycle $\sum \tau_k ^{\prime
\prime}$, which still represents $\phi_*(g_*(\alpha))$ and with the property that {\it the image
of all the singular simplices avoid the point $p$}.

Now observe that all the subdivisions of the singular simplices $\tau_k$ gives rise to 
subdivisions of the singular simplices $\sigma _k$ (recall that we have 
$\tau_k= \phi\circ g\circ \sigma_k$), since a subdivision of a singular simplex 
is actually performed at the level of the
source space. So corresponding to the cycle $\sum \tau_k ^{\prime}$ representing 
$\phi_*(g_*(\alpha))$, we have a corresponding cycle $\sum 
\sigma_k ^{\prime}$ representing the original $\alpha$. Let $\sigma$ be the singular
simplex corresponding to $\tau$, and recall that the subdivision was chosen fine enough so
that $\tau$ was contained inside the image of $g(U)$. In particular, this forces $\sigma (\Delta^n)
\subset U$, which we recall lies in the set of injectivity of the map $g$. This implies that  
there are no ``accidental cancellations'' due to distinct singular simplices in the chain 
$\sum \sigma_k ^{\prime}$ both getting mapped to $\tau$. Since $\tau$ occurred an
even number of times in the cycle $\sum \tau_k ^{\prime}$, we have that $\sigma$ likewise
occurs an even number of times in the cycle $\sum \sigma_k ^{\prime}$. Finally, working with
coefficients in $\mZ_2$ means that we can drop all copies of $\sigma$, obtaining a new cycle
$\sum \sigma_k ^{\prime \prime}$ having the following two properties:
\begin{itemize}
\item the cycle $\sum \sigma_k ^{\prime \prime}$ represents the homology class $\alpha$, and
\item the cycle $\sum \sigma_k ^{\prime \prime}$ has image in $X$ which is disjoint from the point 
$O\in U$.
\end{itemize}
But now observe that, since $U$ is homeomorphic to $\mD^n_\circ$, we have that $X-U$ is a
deformation retract of $X-\{O\}$. Applying the deformation retraction to the cycle
$\sum \sigma_k ^{\prime \prime}$, we can now obtain a cycle representing $\alpha$, but whose
image is contained inside $X-U$. This forces $\alpha \in \im \big( H_n(X-U)\rightarrow 
H_n(X)\big)$, contradicting property (2) in our choice of the class $\alpha$, and completing
our proof of Theorem 2.
\end{Prf}

\section{Optimality and examples}

Before discussing some specific examples covered by our two theorems, let us start by giving some simple non-examples:
\begin{itemize}
\item Take $X=S^1 \subset \mR^2$, and let $f$ be the projection onto the interval 
$[-1,1] \subset \mR ^2 \subset S^2$. The set of injectivity does not contain any open set $U$.
\item Take $X\subset \mR^2$ to be the union of the standard unit circle, along with the interval
$[1,2]$ on the $x$-axis. If $f$ is the projection onto the interval $[-1,2]\subset \mR^2 \subset S^2$,
then we see that $f$ is injective on $U = (1,2)$, but the inclusion $X-U \hookrightarrow X$ induces an
isomorphism on $\check H^1$, so has trivial kernel.
\end{itemize}
In both cases, we see that the conclusion to Theorem 1 fails, i.e. $f(X) \subset S^2$ fails to separate.
Similarly, for Theorem 2, we can consider the following simple example: let $X=S^1$, and 
$f:S^1\hookrightarrow S^2$ be the embedding into the equator. Let $\hat X = S^1 \times [0,1]$ be 
an annulus,
with $X \subset \hat X$ corresponding to $S^1\times \{0\}$. Note that since $X\hookrightarrow \hat X$
is a homotopy equivalence, the group $K$ in Theorem 2 is automatically trivial, and hence 
$K\subset J$ holds. It is also clear that there exist extensions $F: S^1\times [0,1] \rightarrow
S^2$ with the property that $F$ fails to surject onto either hemisphere. These simple examples
show that the hypotheses in Theorems 1 and 2 are indeed necessary.

\vskip 10pt

We now proceed to give some examples of triples $(\hat X, X, U)$ which satisfy the homological 
conditions of both theorems. In particular, for each of the following examples of triples, we have 
that if $f: X\rightarrow S^{n+1}$ is injective on $U$, then: (1) $f(X)$ automatically separates $S^{n+1}$,
(2) there are exactly two connected components $V_1,V_2 \subset S^{n+1}-f(X)$ whose closure
contain $f(U)$, and (3) any extension of $f$ to a map $F:\hat X \rightarrow S^{n+1}$ surjects onto
either $V_1$ or $V_2$.

\vskip 5pt
 
\noindent {\bf Example: manifolds pairs.} 

\vskip 5pt

Taking $\hat X$ to be a compact $(n+1)$-dimensional manifold with non-empty boundary, let
$X$ be the boundary of $\hat X$, and let $U$ be any open $n$-disk in $X$. Note that both of
the groups $\check H^n(X)$ and $H_n(X)$ are free $\mZ_2$-modules, generated by the 
connected components of $X$ (each of which is a closed manifold). It is now immediate 
that the map $\check H^n(X) \rightarrow \check H^n(X-U)$ has non-trivial kernel, as removing
$U$ changes one of the connected components of $X$ from a closed manifold to a compact
manifold with boundary. Similarly, this forces $J=\im\big(H_n(X-U)\rightarrow H_n(X)\big)$ 
to miss elements in $H_n(X)$: an explicit element $\alpha\notin J$ is given by the sum of the 
fundamental classes of the connected components of $X$. Finally, we note that $\alpha \in K$
is the image of the relative fundamental class $\mu \in H_{n+1}(\hat X, X)$
of the manifold with boundary
$\hat X$. From the long-exact sequence of the pair $(\hat X,X)$, we immediately see that
$\alpha$ maps to zero in $H_n(\hat X)$, and hence $\alpha \in K$. This implies that $K\not 
\subseteq J$, completing the verification of the homological hypotheses in our two theorems.

For a more concrete example, if $X=K^2$
denotes the Klein bottle, and if $f: K^2 \rightarrow \mR^3$ is a continuous map which contains
an open set $U\subset Inj(f)$, then the first theorem implies that
$f(K^2)$ separates $\mR^3$ into at least two connected
components. Our second theorem tells us that, if $\hat X= L^3$ is the ``solid Klein-bottle'', and 
$F:L^3\rightarrow \mR^3$ is an extension of $f$, then $F$ surjects onto one of the two components
$V_1,V_2$ incident to any prescribed open set $U\subset Inj(f)$.

The reader might also like to compare the example where $\hat X = S^1 \times [0,1]$ with
the non-example discussed at the beginning of this section. The distinction lies of course in
the choice of $X$, which in the present example is $\partial \hat X = S^1 \times \{0,1\}$, and in 
the non-example, consisted solely of $S^1\times \{0\}$.

\vskip 5pt

\noindent {\bf Example: pseudo-manifold pairs.}

\vskip 5pt

Recall that an $n$-dimensional {\it pseudo-manifold} is a simplicial complex with the property
that every $(n-1)$-dimensional simplex is a face of {\it exactly two} $n$-simplices. 
An $(n+1)$-dimensional pseudo-manifold with boundary is defined to be a simplical
pair $(\hat X,X)$ with
the property that: (1) $X$ is an $n$-dimensional pseudo-manifold, (2) every $n$-simplex
in $X$ is contained in a unique $(n+1)$-dimensional simplex, and (3) every $n$-simplex 
in $\hat X - X$ is contained in exactly two $(n+1)$-simplices. From the homological viewpoint,
the important observation is that the constraint on the codimension one simplices 
ensures that compact pseudo-manifolds still have a notion of a fundamental class (the sum of all
top-dimensional simplices). In particular, the arguments given earlier for manifolds easily
extend to the pseudo-manifold case.

The distinction with the manifold situation is that pseudo-manifolds are allowed to be 
singular, but that the singularities are relatively ``small'', i.e. codimension at least two. The importance
of this class of topological spaces comes from the fact that every complex projective algebraic variety
is a pseudo-manifold (as the singularities will have complex codimension $\geq 1$). 
Note that complex projective algebraic varieties can always be triangulated (see \cite{LW}). Examples
of compact pseudo-manifolds with boundary can be obtained by taking a complex projective variety
$V\subset \mC P^n$, taking a suitable real codimension one, smooth submanifold 
$M^{2n-1} \subset \mC P^n$ which intersects $V$ non-trivially, and cutting $\mC P^n$ open along $M^{2n-1}$. The portion of $V$ in any of the (at most two) connected components of the resulting 
manifold with boundary will yield an example of a compact pseudo-manifold with boundary.

\vskip 5pt

\noindent {\bf Example: CW-complexes and universality.}

\vskip 5pt

We now proceed to consider the case where both $X$ and $\hat X$ are CW-complexes, 
and where $U \subset X$ is a top-dimensional open cell $e ^n \subset X$. In this case, we
note that \v Cech cohomology coincides with singular cohomology.
Furthermore, it is easy to see from a Mayer-Vietoris sequence argument
that the map $\check H^n(X) \rightarrow \check H^n(X-e^n)$ has non-trivial kernel precisely if 
the attaching map $S ^{n-1} \hookrightarrow X^{(n-1)}$ for the $n$-cell $e^n$ induces the zero 
map on $(n-1)$-dimensional cohomology (where $X^{(n-1)}\subset X$ is the $(n-1)$-skeleton
of $X$). As such it is easy to construct CW-complexes satisfying 
the homological conditions of our first theorem. Similarly, it is easy to extend such an 
$n$-dimensional CW-complex $X$ to an 
$(n+1)$-dimensional CW-complex $\hat X$ satisfying the homological conditions for
our Theorem 2. For example, one extension which always works is the case where $\hat X$
is taken to be the cone over $X$. To see this, first observe that $\hat X$ is contractible,
and hence that $H_n(\hat X)=0$. This forces $K= \ker \big(H_n(X) \rightarrow H_n(\hat X)\big)=
H_n(X)$. So as long as the map $H_n(X-u) \rightarrow H_n(X)$ is {\it not} surjective, the
pair $(\hat X, X)$ will satisfy the homological conditions of our Theorem 2. But observe that,
from the fact that $X$ is a CW-complex satisfying the conditions of our Theorem 1,
we have that $H^{n-1}(X-U) \rightarrow 0 \in H^{n-1}(S^{n-1})$.
Since we are working with coefficients in $\mZ_2$, the Ext term in the universal coefficient theorem
automatically vanishes, and we can identify the cohomology groups above as the duals of the
corresponding homology group. This forces the cohomological
statement above to be equivalent to the dual homological statement: 
$\mZ_2 \cong H_{n-1}(S^{n-1}) \rightarrow 0 \in H_{n-1}(X-U)$. The Mayer-Vietoris
sequence now yields:
$$H_n(X-u) \rightarrow H_n(X) \rightarrow \mZ_2 \rightarrow 0$$
confirming that the first map is {\it not} surjective. This completes the verification that taking
$\hat X$ to be the cone over $X$ always satisfies the homological conditions for Theorem 2.

\vskip 10pt

Let us conclude with an observation: the example of pseudo-manifolds
is, in some sense, a ``universal example'' amongst CW-complexes. 
Indeed, let us illustrate what we mean by reconsidering the
situation where $(\hat X, X)$ are CW-complexes, and $U \subset X$ is the interior of an $n$-cell
in $X$.  Since $X$ is assumed to satisfy the hypotheses of Theorem 1, the analysis in the
previous paragraph gives rise to a homology class $\alpha \in H_n(X)$ which is {\it not} in the image
of $H_n(X-U)$. Now recall that, given a homology class $\alpha\in H_n(X; \mZ_2)$, 
in an arbitrary topological space $X$, there exists an $n$-dimensional 
pseudo-manifold $Y$ and a continuous map $\phi: Y\rightarrow X$ with the property that
$\phi_*[Y] = \alpha$.


Now it is easy to see, from the constraints on $\alpha$, that the image $\phi (Y)$ must pass through 
$U$. One can further modify $Y$ so as to ensure that the corresponding
$\phi: Y \rightarrow X$ has the property that $\phi: \phi ^{-1}(U) \rightarrow U$ is a 
homeomorphism (this can be done by an argument similar to the one at the end of Section 2). 
This discussion establishes the following:

\begin{lemma}[Universality of pseudo-manifold example] Let $X$ be an $n$-dimensional  
CW-complex, $U$ an $n$-cell in $X$, such that the pair $(X,U)$ satisfies the hypotheses 
of Theorem 1 for the map $f: X\rightarrow S^{n+1}$. Then there exists a pseudo-manifold 
pair $(Y, U^\prime)$ and a map $\phi: (Y, U^\prime) \rightarrow (X,U)$, with the property 
that the pair $(Y, U^\prime)$ satisfies the hypotheses of Theorem 1 for 
the composite map $f\circ \phi : Y \rightarrow S^{n+1}$.
\end{lemma}

In particular, the image of $\phi(Y)$ under $f$ already separates $S^{n+1}$, showing that
Theorem 1 for CW-complexes is actually a consequence of Theorem 1 for pseudo-manifolds. 
A similar analysis can be used to show that  Theorem 2 for CW-complex pairs can also be 
deduced from the pseudo-manifold case; we leave the details to the interested reader.

\vskip 10pt 

\noindent {\bf Example: spaces which are not CW-complexes.}

\vskip 5pt

Finally, we give an example which {\it cannot} be deduced from the pseudo-manifold case: take
$X$ to be the closed topologists sine curve. It is well known that the first singular cohomology group
is $H^1(X) = 0$, while the first \v Cech cohomology group is $\check H^1(X)\cong \mZ _2$.
If $U\subset X$ is an open interval in the ``sine portion'' of $X$, then one can readily verify that
both $H^1(X-U) = 0$ and $\check H^1(X-U) = 0$. In particular, we see that the pair $(X,U)$ 
satisfy the homological conditions for our Theorem 1, so any continuous map $f: X\rightarrow S^2$
which is injective on $U$ will have image that separates. Of course, the fact that 
$H_1(X) \cong H^1(X) = 0$ tells us that there is no chance of using pseudo-manifolds to detect
separation (in contrast to the situation with CW-complexes).
For $\hat X$, one can again take the cone over $X$; it is easy to verify that the homological
conditions for Theorem 2 are indeed satisfied.

\end{document}